\documentclass[11pt,twoside]{amsart}
\newcounter{alphthm}
\newtheorem{propriete}[alphthm]{Theorem}

\def\beq{\begin{equation}}
\def\eeq{\end{equation}}
\def\nn{\nonumber}

\newtheorem{thm}{Theorem}[section]
\newtheorem{cor}[thm]{Corollary}

\newtheorem{defn}[thm]{Definition}
\newtheorem{rem}[thm]{\bf{Remark}}

\numberwithin{equation}{section}
\def\pn{\par\noindent}

\begin{document}

\leftline{ \scriptsize \it Bulletin of the Iranian Mathematical
Society  Vol. {\bf\rm XX} No. X {\rm(}200X{\rm)}, pp XX-XX.}

\vspace{1.3 cm}

\title
{Properties of  Generalized Berwald Connections}
\author{ B. Bidabad* and  A. Tayebi}
\thanks{{\scriptsize
\\
\hskip -0.4 true cm MSC(2000): Primary: 53B40; Secondary: 53C60.
Received:, Accepted:\\
$*$Corresponding author
\newline\indent{\scriptsize $\copyright$ 2008 Iranian Mathematical
Society}}}
\maketitle
\begin{center}
\end{center}
\begin{abstract}  Recently the present authors introduced a general class of Finsler
connections which leads to a smart representation of connection
theory in Finsler geometry and yields to a classification of Finsler
connections into the three classes. Here the properties of one of
these classes namely the Berwald-type connections which contains
Berwald and Chern(Rund) connections as a special case is studied. It
is proved among the other that the hv-curvature of these connections
vanishes if and only if the Finsler space is a Berwald one.  Some
applications of this connection is discussed.
\end{abstract}
\vskip 0.2 true cm
\pagestyle{myheadings}
\markboth{\rightline {\scriptsize  B. Bidabad and  A. Tayebi}}
         {\leftline{\scriptsize Properties of  Generalized Berwald Connections}}
\bigskip
\bigskip
\vskip 0.4 true cm
\section{\bf Introduction}
\vskip 0.4 true cm
\setcounter{equation}{0}
 Always there were a hope to find a solution
for some of the unsolved problems by developing a connection theory,
thus it is  useful  to introduce new connections in Finsler
geometry. As it is mentioned in \cite{MS}, the study of hv-curvature
of Finsler connections is of urgent necessity for the Finsler
geometry as well as for theoretical physics. Similarly as another
application of Finsler connections in physics one can mention an
example in Relativistic field theory.  In this theory different
connections have been defined in Finsler geometry where the
connections, torsions, or curvatures can be related to fields which
might be identified as electromagnetic or Yang-Mills fields. In this
relation, one can refer to \cite{Ak2}, \cite{Be} and \cite{MA}.

Historically, in Riemannian geometry, the connection of choice was
that constructed by Levi-Civita,  which has two remarkable
attributes; metric-compatibility and torsion-freeness. In 1926, L.
Berwald \cite{Be} introduced a connection and two curvature tensors.
The Berwald connection is torsion-free, but is
 not necessarily metric-compatible. It was Berwald who first
successfully extended the notion of Riemann curvature to Finsler
spaces. He also introduced a notion of non-Riemannian quantity
called Berwald curvature. From this point of view, Berwald is the
founder of differential geometry of Finsler spaces \cite{Sh2}. Next,
Cartan in 1934 has find locally the coefficients of a
metric-compatible and h-torsion free connection, called later,
Cartan connection. The global construction of this connection is
given in a remarkable work of Akbar-Zadeh in 1967 \cite{Ak1}. Other
progress came in 1948, when the Chern (Rund) connection was defined.
In 1943 Chern studied the equivalence problem for Finsler spaces
using the Cartan exterior differentiation method \cite{Ch}. Chern
came back to his connection in 1993, in a joint paper with Bao
\cite{BC} and shows its usefulness in treating global problems in
Finsler geometry.
 The Berwald and Chern connections also fail slightly but expectedly,
 to be
metric-compatible. The Chern connection has a simpler form, while
the Berwald connection affects a leaner hh-curvature for spaces of
constant flag curvature. Indeed the Berwald connection is
particularly convenient when dealing with Finsler spaces of constant
flag curvature. It is most directly related to the nonlinear
connection coefficients and most amenable to the study of the
geometry of paths. These connections (Berwald and Chern) coincide
when the underlying Finsler structure is of
 Landsberg type. They further reduce to a linear connection on $M$,
  when the Finsler structure is of Berwald type \cite{CS}, \cite{M}.

Recently the present authors have  defined a  general class of
Finsler connections which leads to a general representation of some
Finsler connections  in Finsler geometry and yields to a
classification of Finsler connections into the three classes,
namely, Berwald-type, Cartan-type and Shen-type connections
\cite{BT}. In the present work we study the properties of the former
connection which
  contains  Berwald and Chern(Rund) connections  as special cases
  and is the most general connection of this kind. We prove in continuation of Berwald's and Chern's  works that the hv-curvature of the Berwald-type connection characterizes the Berwald structure.

A distinguished property of the introduced  connection is its
adaptive form for different applications. In fact one can use a
suitable special case of this connection to find a geometric
interpretation for solutions of certain differential equations
formed by Cartan tensor and its derivatives. For example in section
4 we prove: Let (M,F) be a complete Finsler manifold with bounded
Landsberg tensor. Then $F$ is a Landsberg  metric if and only if its
hv-curvature $P_{jkl}$ vanishes.
\section{\bf {\bf \em{\bf Preliminaries}}}
 Let $M$ be a n-dimensional $ C^\infty$ manifold. For a point $x \in M$,
 denoted by $T_x M $
the tangent space of M at $x$. The tangent bundle of M is the union
of tangent spaces. $TM:=\cup _{x \in M} T_x M$. We will denote
elements in TM by $(x, y)$ if $y\in T_xM$. Let $TM_0 = TM \setminus
\{ 0 \}.$ The natural projection $\pi: TM
\rightarrow  M$ is given by $\pi (x,y):= x$.\\
 Throughout this paper, we use {\it
Einstein summation convention} for expressions with
indices.\footnote{That is where ever an index is appeared twice as a
subscript as well as a superscript, then that term is
assumed to be summed over all values of that index.}\\

A {\it Finsler structure} on a manifold $M$ is a function $ F:TM
\rightarrow [0,\infty )$
with the following properties:\\
(i)$F$ is $C^\infty$ on $TM_0$.\\
(ii)$F$ is positively 1-homogeneous on the fibers of tangent
 bundle $TM$:
 \[
 \forall \lambda>0 \quad  F(x,\lambda y)=\lambda F(x,y).
 \]
 (iii) The Hessian of $F^{2}$ with elements
$ g_{ij}(x,y):=\frac{1}{2}[F^2(x,y)]_{y^iy^j} $ is positively
defined on $TM_0$.\\
Then the pair  $(M,F)$ is called a {\it Finsler manifold.} $F$ is
Riemannian if
$g_{ij} (x,y)$ are independent of $y \neq 0$.\\\\
{\bf{Nonlinear connection.}}\\
Let us consider the tangent bundle $(TM, \pi, M)$ of the manifold
$M$. The tangent bundle of the manifold $TM$ is $(TTM, \pi_*, TM)$,
where $\pi_*$ is the tangent mapping of the projection $\pi$. A
tangent vector field on $TM$ can be represented in the local natural
frame $({{\partial }\over {\partial x^i}}, {{\partial }\over
{\partial y^i}})$ on $TM$ by
\[
\bar X=X^i(x,y){{\partial }\over {\partial x^i}}+Y^i(x,y){{\partial
}\over {\partial y^i}}.
\]
It can be written in the form $\bar X=(x, y, X^i, Y^i)$ or, shorter,
$\bar X=(x, y, X, Y)$. The mapping $\pi_*:TTM \rightarrow TM$ has
the local form
\[
\pi_*(x, y, X, Y)=(x, y).
\]
Put $ VTM:=ker \pi_*=span \{{{\partial }\over {
\partial y^i}} \} ^n_{i=1} $. $VTM$ is an n-dimensional sub-bundle
of $T(TM_0)$, whose fiber $V_v TM$ at {\it v}  is just the tangent
space $T_v (T_x M)\subset T_v (TM_0)$.  $VTM$ is called the {\it
vertical tangent bundle} of $TM_0.$ \\
We can write the vertical subbundle as $(VTM, \pi_{_{VTM}}, TM)$.
Its fibres are the linear vertical spaces $V_uTM$, $u \in TM.$ The
points of submanifold $VTM$ are of the form $(x, y, 0,Y )$. Hence,
the fibers $V_uTM$ of the vertical bundle are isomorphic
to the real vector space $\mathbb{R}^n$.\\

Let us consider the pullback tangent bundle $\pi^*TM$ defined as
follows (\cite{BCS}):
\[
\pi^*TM=\left\{(u, v)\in TM \times TM | \pi(u)\\=\pi(v)\right\}.
\]
Take a local coordinate system $(x^i)$ in M. The local natural frame
$\{{{\partial} \over {\partial x^i}}\}$ for $T_xM$ determines a
local natural frame $\partial_i$ for $\pi^*TM$, ${
\partial _i |_y
:= (y,{{\partial} \over {\partial x^i}}| _x )}$, $y \in T_xM$. This
gives rise to a linear isomorphism between $\pi^*TM|_y$ and $T_xM$
for every $y \in T_xM.$ There is a canonical section $\ell$ of $\pi
^*TM$ defined  by
$\ell=\ell^i\partial_i$, where $\ell^i=y^i/{F(x,y)}$.\\
 The fibers of $\pi^*TM$, i.e., $\pi^*_uTM$ are isomorphic to
 $T_{\pi(u)}M$. One
can define the following morphism of vector bundle $\rho: TTM
\rightarrow \pi^*TM$, $\rho(X_u)=(u, \pi_*(\bar X_u))$. It follows
that
\[
ker \rho=ker \pi_*=VTM
\]
By means of these consideration one can see without any difficulties
that the following sequence is exact
\begin{equation}
0{\rightarrow} VTM \overset{i}{\rightarrow} TTM
\overset{\rho}\rightarrow \pi^*TM \rightarrow 0 \ ,
 \end{equation}
where $i$ is natural inclusion map.\\

A {\it nonlinear connection} on the manifold $TM$ is a left
splitting of the exact sequence (1.1). Therefore, a nonlinear
connection on TM is a vector bundle morphism $C: TTM \rightarrow
VTM$, with the property $C\circ i=1_{_{VTM}}$. The kernel of the
morphism $C$ is a vector bundle of the tangent bundle $(TTM, \pi_*,
TM)$, denoted by $(HTM, \pi_{_{HTM}}, TM)$ and called the {\it
horizontal subbundle}. Its fibres $H_uTM$ determine a distribution
$u \in TM\rightarrow H_uTM\subset T_uTM$, supplementary to the
vertical distribution $u \in TM\rightarrow V_uTM\subset T_uTM$.
Therefore, a nonlinear connection N induces the following Whitney
sum:
\begin{equation}
TTM=HTM \oplus VTM.
\end{equation}
Let we put
 \begin{equation}
\frac{\delta}{\delta x^j}:=\frac{\partial}{\partial
x^j}-N^i_j\frac{\partial}{\partial y^i}\ ,
 \end{equation}
 where the above $N^i_j$ are the components of $N$ and are known
in the trade as the {\it nonlinear
    connection coefficients} on $TM_0$.\\
Restriction of the morphism $\rho: TTM\rightarrow \pi^*TM$ to the
$HTM$ is an isomorphism of vector bundles, for which we have
 \begin{equation}
   \rho ({\partial \over {\partial
x^i}})=\partial_i \quad,\quad \rho({\partial \over {\partial
y^i}})=0.
  \end{equation}
Let $\nabla$ be a linear connection on $\pi^*TM$, $\nabla : \chi(TM_0)\times \pi^*TM \rightarrow \pi^*TM$ such that $\nabla:
(\hat{X},Y) \rightarrow \nabla_{\hat{X}}Y$. A {\it Finsler
connection} is a pair of a linear connection $\nabla$, and a
nonlinear connection $N$.\\

Given a Finsler metric F on M, $F(y)=F(y^i{{\partial} \over
{\partial x^i}}|_{x})$ is a function of $(y^i) \in \mathbb{R}^n$ at
each point $x \in M$.  Finsler metric $F$ defines a fundamental
tensor $g:\pi ^* TM \otimes \pi ^* TM \rightarrow [0,\infty)$ by the
formula $g(\partial_i |_v,\partial_j |_{v})=g_{ij}(x,y)$, where
$v=y^i {{\partial} \over {\partial x^i}}|_x$ and $g_{ij}$
 are defined in the definition of Finsler structure.
Then $(\pi ^* TM,g)$ becomes a Riemannian vector bundle over $TM_0$.
Let
$$
 A_{ijk}(x,y)={1 \over 2}F(x,y)[F^2(x,y)]_{y^iy^jy^k}.
$$
Clearly, $A_{ijk}$ is symmetric with respect to $i,j,k$. The Cartan
tensor $A :\pi ^* TM \otimes \pi ^* TM \otimes\pi ^* TM \rightarrow
\mathbb{R} $ is defined by $ A(\partial_i |_v ,\partial_j |_
{v},\partial_k |_ {v})=A_{ijk}(x,y). $ In some literatures
$C_{ijk}=A_{ijk}/F $ is called Cartan tensor. Riemannian manifolds
are characterized by $A\equiv0.$ F is positively homogenous of
degree 1 on M then by the Euler's theorem we see that $y^iF_{y^i}=F$
and then $y^iF_{y^iy^j}=0$ using this the canonical section $\ell$
satisfies:
\[
g(\ell,\ell)=1 \quad , \quad  A(X,Y,\ell)=0,
\]
where the second equation is equivalent to
$A(X,Y,\frac{y^i}{F}\frac{\partial}{\partial x^i})=\ell^i
A(X,Y,\frac{\partial}{\partial x^i})=0$.  Let $\bar \ell$ denote the
unique vector field in $HTM$ such that $\rho(\bar \ell)=\ell$. We
call $\bar \ell$ a {\it geodesic } or {\it spray} field on $TM_0 $.\\

 Let $\nabla$ be the Berwald (or Chern) connection,
by means of $\nabla$, the tensor  $\dot A$ is defined by $\dot A :
\pi ^* TM \otimes \pi ^* TM \otimes \pi ^* TM \rightarrow
\mathbb{R}$,
\[
\dot A (X, Y, Z) := \bar \ell A(X,Y,Z) - A(\nabla_{\bar \ell} X, Y,
Z)-
 A(X,\nabla_{\bar \ell} Y, Z)- A(X,Y,\nabla_{\bar \ell} Z).
 \]
 Putting $\overset{_1}{A}_{ijk}=\dot A_{ijk}$, $\overset{_2}{A}_{ijk}=\ddot
 A_{ijk}$,
 $\forall m\in \mathbb{N}$ we define
 $\overset{_{m+1}}{A}$ as follow
\[
 \overset{_{m+1}}{A}(X, Y, Z) := \bar \ell \quad \overset{_m}{A}(X,Y,Z)
  -\overset{_m}{A}
 (\nabla_{\bar \ell} X,
Y, Z)- \overset{_m}{A} (X,\nabla_{\bar \ell} Y, Z)
 - \overset{_m}{A}
(X,Y,\nabla_{\bar \ell} Z).\quad \quad \quad \quad \quad \quad \quad
\quad \quad
\]
Obviously, $\forall m\in\mathbb{N}$, the tensors
$\overset{_m}{A}_{ijk}$ are symmetric with respect to three indices.
Moreover, using $\nabla_{\bar \ell}\ \ell=0 $ we have
$\overset{_m}{A}(X,Y,\ell)=0$, $\forall m\in\mathbb{N}$. $A$ and
$\dot A$  are basic tensors in Finsler
 geometry. In the Riemannian case, both of them vanish. Therefore by the
 above definition we know that in the Riemannian case $\forall
 m\in\mathbb{N}$,  $\overset{_m}{A}=0.$\\

 A Finsler metric $F(x,y)$ on a manifold $M$ is called {\it Berwald
metric} if in any standard local coordinate system $(x^i,y^i)$ in
$TM_0$, the Christoffel symbols $\Gamma^k _{ij}=\Gamma^k _{ij}(x)$
are functions of $x \in M$ alone. In this case
$G^i(x,y)=\frac{1}{2}\Gamma^i _{jk}(x)y^jy^k$ are quadratic in
$y=y^i {{\partial} \over {\partial x^i}}|_x$ and $F(x,y)$ is called
a {\it  Landsberg  metric} if $L^i_{jk}(x,y)=0$, that is
$$
L^i_{jk}(x,y)=\frac {\partial^2 G^i}{\partial y^j \partial
y^k}(x,y)-\Gamma^i _{jk}(x,y).
$$
Clearly Minkowski and Riemannian metrics are trivial Berwald
metrics. If $F(x,y)$ is a Berwald metric, it is a Landsberg metric.
But the converse might not be true, although no counter-example has
been found yet \cite{Sh2}. A fundamental theorem in Finsler geometry
says that a Finsler metric $F$ is a Berwald metric if and only if
the Cartan tensor is covariantly constant along all horizontal
directions on the slit tangent bundle $TM_0$ (see \cite{Sh3} and
\cite{BCS} for proof). Thus in the Berwald case, $
\overset{_m}{A}_{ijk}$
vanish $\forall m\in \mathbb{N}$.\\\\
{\bf Flag curvature}. A flag curvature is a geometrical invariant
that generalizes the sectional curvature of Riemannian geometry. Let
$x \in M$, $0\neq y \in T_xM$ and $V:=V^i \frac{\partial}{\partial
x^i}$. Flag curvature is obtained by carrying out the following
computation at the point $(x, y)\in TM_0$, and viewing y, V as
section of $\pi^*TM$:
\[
K(y, V):=\frac{V^i(y^j \ R_{jikl}\ y^l) V^k}{g(y,y)g(V,V)-[g(y,
V)]^2} ,\] where $g$ is a Riemannian metric on $\pi^*TM$. If $K$ is
independent of the transverse edge $V$, we say that our Finsler
space has {\it scalar flag curvature}. Denote this scalar by
$\lambda=\lambda(x, y)$. When $\lambda(x, y)$ has no dependence on
either $x$ or $y$, then Finsler manifold is said to be of {\it
constant flag curvature.}

\vskip 0.4 true cm

\section{\bf {\bf \em{\bf Berwald-type connection on
$\pi ^* TM$}}}

\vskip 0.4 true cm

In this section we introduce
a new family of Finsler connections which are torsion-free and
almost compatible with the Finsler metric. In the sequel we will
refer to this connection by {\it ``Berwald-type connection"}.
\begin{defn}
\emph{Let $(M,F)$ be a Finsler n-manifold. Let $g$ and
 $A$ denote the fundamental and the Cartan tensors in $\pi ^* TM$,
respectively. Let $D$ be a Finsler connection on $M$.}\\
 (i) $D$ is {\bf torsion-free}, if $
\forall \hat{X},\hat{Y} \in {\chi X} (TM_0),$
\begin{equation}{ \mathfrak{T}} _D (\hat{X},
\hat{Y}):=D_{\hat{X}} \rho(\hat{Y})-D_{\hat{Y}}
  \rho(\hat{X})-\rho([\hat{X},\hat{Y}])=0.
\end{equation}
\emph{(ii) $D$ is {\bf almost compatible} with the Finsler structure
in the following
  sense: if for all} $X,Y \in
\pi ^* TM$ and $\hat{Z} \in T_v (TM_0)$,
\begin{eqnarray}
\nonumber(D_{\hat{Z}} g)(X,Y):= \hat{Z}g(X,Y) -g(D_{\hat{Z}} X, Y)
-g(X, D_{\hat{Z}} Y),
 \end{eqnarray}
\emph{or equivalently}
\begin{eqnarray}
\nonumber(D_{\hat{Z}} g)(X,Y)&=&-2k_{_1}
 \dot A(\rho(\hat{Z}),X,Y)-\cdots-2k_{_m}
 \overset{_m}{A}(\rho(\hat{Z}),X,Y)\\&+&2F^{-1}
A(\mu(\hat{Z}),X,Y),
 \end{eqnarray}
\emph{where $\rho(\hat Z ):=(v,\pi_* (\hat Z))$, $\mu (\hat
Z):=D_{\hat Z} F  \ell$, $m \in \mathbb{N}$ and $k_{_{i}}\in
\mathbb{R}$.}
\end{defn}
\vskip 0.4 true cm

The bundle map $\mu:T(TM_0) \rightarrow \pi ^* TM$ defined in above
definition satisfies
 $
   \mu({\partial \over \partial y^i})=\partial_i.
   $
To prove this, take $\hat \ell=\ell^i{ \partial \over
\partial x^i}$, where $\ell=\ell ^i \partial_i$. Now $ \rho(\hat
\ell)=\ell$, so from (2 .1)
  \begin{equation}
   \mu({\partial \over \partial y^i})=D_{\partial \over \partial
y^i} F \ell = \rho ([{\partial \over \partial y^i},y^k {\partial
\over \partial x^k}])=\partial_i.
   \end{equation}
\setcounter{thm}{0}\begin{thm} Let $(M,F)$ be a Finsler n-manifold.
Then there is a unique linear torsion-free connection $D$ in $\pi ^*
TM$, which is almost compatible with the Finsler structure in the
sense of (2 .2).
\end{thm}
\vskip 0.4 true cm

\begin{proof} In a standard local coordinate system $(x^i , y^i)$
in $TM_0$, we write
\[
D_{{{\partial} \over {\partial x^i}}} {\partial_j} =\Gamma^k _{ij}
{\partial_k}\quad, \quad D_{{{\partial} \over {\partial y^i}}}
{\partial_j} =F^k _{ij} {\partial_k}\ .
\]
By replacing $\hat X, \hat Y$ in (2 .1) with the basis of
$T_v(TM_0)$ i.e. $\{\frac{\partial}{\partial
x^i},\frac{\partial}{\partial y^i}\}$  we get
\begin{equation}
 \Gamma ^k _{ij}=\Gamma ^k _{ji},\quad \quad \quad \quad
   \quad\\
\end{equation}
\begin{equation}
  F ^k _{ij}=0, \quad \quad \quad \quad \quad \quad\\
\end{equation}
and by replacing $X, Y$ (resp. $\hat Z$ ) in (2 .2) with the basis
of $\pi^*TM$ i.e. $\{\partial_i\}$ (resp. with the basis of
$T_v(TM_0)$) we get
  \begin{eqnarray}
\nonumber{{\partial}\over {\partial x^k}} (g_{ij})={\Gamma^l} _{ki}
g_{lj}+\Gamma ^l _{kj}g_{li}&-&2k_{_1}\dot
A_{ijk}-\cdots-2k_{_m}\overset{_m}{A}_{ijk}
\\  &+&2A_{ijl}\Gamma ^l
_{km} \ell^m,\\
\nonumber {{\partial}\over {\partial y^k}} (g_{ij})={F^l} _{kj}
g_{li}+F^l _{ik}g_{jl}&-&2\{k_{_1}\dot A_{ijk}+\cdots
+k_{_m}\overset{_m}{A}_{ijk}\}F^l _{mk} \ell^m\\ &+&2 F^{-1}A_{ijk},
\end{eqnarray}
 where
$g_{ij}, A_{ijk}$ and $A^{^{(m)}}_{ijk}$, $\forall m\in\mathbb{N}$
are all functions of $(x,y)$. We shall compute $\Gamma ^k _{ij}$ by
" Christoffel's trick " from (2 .4) and (2 .6). Then making a
permutation to $i, j, k$ in (2 .6), and using (2 .4), we obtain
\begin{eqnarray}
 \nonumber \Gamma ^k _{ij}&=&\gamma ^k _{ij} +k_{_1}\dot A^k_{\ ij}
+\cdots+k_{_m} \overset{_m}{A ^k}_{ij}\\ &+&g^{kl} \left\{ A _{ijm}
\Gamma ^m _{lb}-A _{jlm} \Gamma ^m _{ib}-A _{lim} \Gamma ^m
_{jb}\right\} \ell^b,
  \end{eqnarray}
where we have put
 \[
   \gamma ^k _{ij}={1 \over 2}g^{kl}\left\{{{\partial g_{jl}}\over
   {\partial
x^i}}
   + {{\partial g_{il}}\over {\partial x^j}} -
    {{\partial g_{ij}}\over {\partial x^l}} \right\},
     \]
     and $A^k _{ij}=g^{kl}A _{ijl}$. Multiplying (2 .8) by $\ell^i$,
we obtain
 \begin{equation}
 \Gamma ^k _{ib}\ell^b = \gamma ^k _{ib}\ell^b -A^k _{\ im} \Gamma
^m
 _{lb}\ell^l \ell^b.
 \end{equation}
 Multiplying (2 .9) by $\ell^j$, yields
 \begin{equation}
 \Gamma ^k _{ab} \ell^a \ell^b = \gamma ^k _{ab} \ell^a \ell^b.
 \end{equation}
 Substituting (2 .10) into (2 .9), we obtain
 \begin{equation}
         \Gamma ^k _{ib}\ell^b = \gamma ^k _{ib} \ell^b - A^k
_{\ im} \gamma ^m
 _{ab}\ell^a \ell^b.
 \end{equation}
 Substituting (2 .11) in (2 .8), we obtain
\begin{eqnarray}
\Gamma ^k _{ij} =\gamma ^k _{ij}&+&k_{_1}\dot A^k_{\ ij}
+\cdots+k_{_m} \overset{_m}{A ^k}_{ij} +g^{kl} \left\{ A _{ijm}
\gamma ^m _{lb}-A _{jlm} \gamma ^m _{ib}-A
_{lim} \gamma ^m _{jb}\right\} \ell^b \nn \\
&+&\left\{ A^k _{jm}A^m _{is}+A^k _{im}A^m _{js}-A^k _{sm}A^m
_{ij}\right\} \gamma ^s _{ab}\ell ^b \ell^a.
 \end{eqnarray}
 Then using (1 .3), (2 .12) become
\begin{equation}
\Gamma ^i _{jk} =\frac{g^{is}}{2}\{\frac{\delta g_{sj}}{\delta
x^k}-\frac{\delta g_{jk}}{\delta x^s}+\frac{\delta g_{ks}}{\delta
x^j}\}+k_{_1}\dot A^i_{\ jk} +\cdots+k_{_m} \overset{_m}{A ^i}_{\
jk}.
\end{equation}
 This proves the uniqueness of $D$. The set $\{\Gamma ^k _{ij} , F^k
 _{ij}=0\}$ where $\{\Gamma ^k _{ij}\}$ are given by (2 .13), define
  a linear connection $D$ on $\pi^*TM$ satisfying (2 .1) and
 (2 .2).
  \end{proof}

\begin{defn} Let $(M,F)$ be a Finsler manifold. A Finsler
connection is called of {\it Berwald-type} (resp. {\it of
Cartan-type} or {\it Shen-type}) if and only if  vanishing of its
hv-curvature,  reduces the Finsler structure  to the Berwaldian
(resp. Landsbergian or Riemannian) one.
\end{defn}
\vskip 0.4 true cm

From this  view point one can compare some of the non-Riemannian
Finsler connections according to the compatibility of the tensors
$S$ and $T$.
\begin{center}  \textbf{ A classification of Finsler
connections  according  to their compatible tensors $S$ and
$T$}\end{center} {\small
 \centerline
{\begin{tabular} {||l||c||c||c||c||} \hline\hline
&\multicolumn{2}{|c||}{\textbf{Compatible tensors}}&\ &
\\ \hline \textbf{Connection}& \textbf{S}&\textbf{T}& \textbf{Metric compatibility} & \textbf{Torsion}\\
 \hline \hline\hline 1. Berwald & $A+\overset{\bullet}{A}$ &
  \quad $0$& almost compatible&free\\
\hline 2. Chern- Rund &$A$ &\quad$0$&  almost compatible& free \\
\hline  3. Berwald-type & $A+\kappa_{_1}
\overset{\bullet}{A}+\cdots+\kappa_{_m}\overset{_m}{A}$ &\quad $0$&
almost
 compatible&  free\\
\hline \hline \hline 4.  Cartan & $A$ & \quad$A$& metric compatible&
not free  \\
\hline  5.  Hashiguchi &$A+\overset{\bullet}{A}$ & \quad$A$& almost
compatible& not free\\
\hline  6. Cartan-type &$A+\kappa_{_1}
\overset{\bullet}{A}+\cdots+\kappa_{_m}\overset{_m}{A}$  \ \ \ \ \ \
\  &\ \ \
$A$& depends on $\kappa_{_i}$& not free\\
\hline \hline\hline  7. Shen & $0$ & \quad$0$& almost compatible& free\\
\hline  8. Shen-type& $\kappa_{_1}
\overset{\bullet}{A}+\cdots+\kappa_{_m}\overset{_m}{A}$ & \quad$0$&
almost compatible& free\\
\hline \hline \hline 9. General-type  &$\kappa_{_0} A+\kappa_{_1}
\overset{\bullet}{A}+\cdots+\kappa_{_m}\overset{_m}{A}$ &\quad$r A$&
depends on
$\kappa_{_i}$ and $r$& \footnotesize{depends on $r$}\\
\hline
\end{tabular}}}
\bigskip
 In this table $A$, $\dot A$, $\ddot A$,...,
$\overset{_m}{A}$ are Cartan tensor and their covariant derivatives,
$\kappa_{_i}$ and $r$ are arbitrary real constants. The connections
1, 2, and 3 belong to the Berwald-type category. The connections 4,
5, and 6 are Cartan-type connections.
 The connections 7 and 8 belong to the Shen-type Category. The
 connection 9 contains all other connections.

 \setcounter{thm}{0}
 \begin{rem}\emph{ The Berwald and Chern connections are
 special  cases of Berwald-type connection in the following way:\\
 Putting $k_{_{1}}=\cdots=k_{_{m}}=0$  yields the
 {\bf Chern connection}.\\
Putting $k_{_{2}}=\cdots=k_{_{m}}=0$ and $k_{_{1}}=1$ yields the
{\bf Berwald connection}.}
\end{rem}

 The bundle map $\mu : T(TM_0) \rightarrow \pi ^* TM$  defined in
Definition 1
 can be expressed in the following form;
\begin{equation}
   \mu ({{\partial} \over {\partial x^i}})=N^k _i \partial_k,
   \quad \quad \mu ({{\partial} \over {\partial y^i}})=\partial_i,
\end{equation}
    where $N^k _i =F \Gamma ^k _{ij} \ell ^j=F\{\gamma ^k _{ij} \ell ^j
-A^k _{\ il}
    \gamma ^l _{ab}\ell^a \ell ^b\}$.\\
Using the nonlinear connection coefficients, for Berwald-type
connection we have
\begin{equation}
\Gamma ^i _{jk}=\gamma^i_{jk}+k_{_1}\dot A^i_{\ jk} +\cdots+k_{_m}
\overset{_m}{A ^i}_{\
jk}-g^{il}\{C_{ljs}{N^s_k}-C_{jks}{N^s_l}+C_{kls}{N^s_j}\}.
\end{equation}
We summarize that $Ker \rho =VTM$, $Ker \mu =HTM$, $\rho$ restricted
to $HTM$ is an isomorphism onto $\pi ^* TM$, and $\mu$ restricted to
$VTM$ is the bundle isomorphism onto $\pi ^* TM.$

\section{\bf {\bf \em{\bf Curvature tensors}}}
\setcounter{footnote}{0} In this section we study the curvature
tensors of the Berwald-type connection. This connection is
torsion-free and almost compatible with Finsler metric in the sense
of (2 .2). As a torsion-free connection, it defines two curvatures
$R$ and $P$. The $R$-term is the so-called Riemannian curvature
tensor which is a natural extension of the usual Riemannian
curvature tensor of Riemannian metrics, while the $P$-term is purely
non-Riemannian quantity. We prove also that the hv-curvature $P$ of
this connection vanishes if and only if the Finsler structure is a
Berwald structure.
 The curvature tensor $\Omega$ of $D$ is
defined by
 \begin{equation}
  \Omega(\hat X,\hat Y)Z=D_{\hat X}{D_{\hat Y}Z}-D_{\hat Y}{D_{\hat
X}Z}-
  D_{[\hat X,\hat Y]} Z ,
  \end{equation}
where $\hat X,\hat Y \in \chi (TM_0)$ and $Z \in \pi^* TM.$\\

Let $\{e_i\}^n _{i=1}$ be a  local orthonormal (with respect to $g$)
frame field for the vector bundle $\pi ^* TM$ such
  that $g(e_i ,e_n )=0, i=1,...,n-1$ and $e_n :=\ell$.
Put $\ell_i:=g_{ij}\ell^j=F_{y^i}$. Let $\{\omega^i\}^n_{i=1}$ be
its dual co-frame field. The $\omega^i$s are local sections of the
dual bundle $\pi^* TM$. One readily finds that $
\omega^n:{\partial{F} \over {\partial
  {y^i}}}dx^i=\omega$, which is the {\it Hilbert form}. It is obvious that
  $\omega(\ell)=1$.

   Put
  \[
  \rho =\omega^i \otimes e_i,\quad  De_i = \omega ^{\ j} _i \otimes
e_j, \quad
  \Omega e_i=2\Omega ^{\ j} _i \otimes e_j.
  \]
  $\{\Omega ^{\ j} _{i}\}$ and $\{\omega ^{\ j} _{i}\}$ are called the
  {\it curvature forms} and {\it connection forms} of $D$ with respect to
  $\{e_{i}\}$.
  We have $\mu:=DF\ell=F\{\omega ^{\ i} _n +d(log F)\delta ^i _n
\}\otimes
  e_i.$ Put $\omega^{n+i}:=\omega^{\ i} _n  +d(log F)\delta ^i _n.$ It
is easy to see that $\{\omega ^i, \omega^{n+i} \}^n _{i=1}$ is a
  local basis for $T^*( TM_0).$ By definition
\[
   \rho =\omega ^i \otimes
  e_i,\quad \mu=F \omega^{n+i}\otimes
  e_i.
\]
According to Theorem 1 there exits a
 connection 1-forms $\{\omega^i_j\}$ which satisfy the following
torsion-freeness
 and almost compatibility as
follows.
 \begin{equation}
   d\omega ^i =\omega ^j \wedge \omega^{\ i} _j,
 \end{equation}
   \begin{equation}
dg_{ij}=g_{kj} \omega ^{\ k} _i +g_{ik } \omega ^{\ k} _j
-2\{k_{_1}\dot A_{ijk}+\cdots+k_{_m}\overset{_m}{A}_{ijk}\} \omega
^k +2A_{ijk} \omega ^{n+k}.
  \end{equation}
In fact using the local orthonormal frame field $\{e_i\}^n _{i=1}$
for the vector bundle $\pi ^* TM$ and
 its dual co-frame field $\{\omega^i\}^n_{i=1}$, (2 .1) and
(2 .2) respectively, after a straightforward calculation
analogous to the proof of Theorem 1, become (3 .2) and (3 .3).\\
Let us we put
  \begin{equation}
dg_{ij}-g_{kj} \omega ^{\ k} _i -g_{ik } \omega ^{\ k}
_j=g_{ij|k}\omega^k +g_{ij.k}\omega^{n+k},
    \end{equation}
 where  $g_{ij.k}$ and $g_{ij|k}$ are respectively the vertical
 and horizontal covariant derivative of $g_{ij}$. This gives
 \begin{equation}
g_{ij|k}=-2\{k_{_1}\dot
A_{ijk}+\cdots+k_{_m}\overset{_m}{A}_{ijk}\},
 \end{equation}
and
   \begin{equation}
g_{ij.k}=2 A_{ijk}.
  \end{equation}
 Moreover the torsion freeness is equivalent to following
 \begin{equation}
  \omega ^{\ i} _j=\Gamma^i_{jk}dx^k.
  \end{equation}
 Clearly (3 .1) is equivalent to
   \begin{equation}
   d\omega ^{\ j} _i -\omega ^{\ k}_i \wedge \omega^{\ j}_ k=\Omega^{\
j}_ i.
   \end{equation}
Since the $\Omega ^{\ i} _j$ are 2-forms on the manifold $TM_0$,
they can be generally expanded as
  \begin{equation}
\Omega^{\ j}_ i={1 \over 2}R^{\ j} _{i \ kl} \omega ^k \wedge
\omega^l +P^{\ j} _{i \ kl} \omega ^k \wedge \omega^{n+l}+{1 \over
2}Q^{\ j} _{i \ kl} \omega ^{n+k} \wedge \omega^{n+l},
  \end{equation}
  Let $\{\bar e_i, \dot e_i\}^n _{i=1}$ be the local basis
  for $T(TM_0)$, which is dual to $\{\omega ^i, \omega^{n+i} \}^n
  _{i=1}$, i.e., $\bar e_i \in HTM, \dot e_i \in VTM$ such that
  $\rho(\bar e_i)=e_i, \mu(\dot e_i)=F e_i$.
The objects $R$, $P$ and  $Q$ are respectively the hh-, hv- and
vv-curvature
   tensors of the connection $D$ and with $R(\bar e_k,\bar
e_l)e_i =R^{\ j}
   _{i \ kl}e_j , \quad P(\bar e_k,\dot e_l)e_i=P^{\ j}_{i \ kl}
   e_j,
\quad$
   and $ \quad Q(\dot e_k,\dot e_l)e_i=Q^{\ j} _{i \ kl} e_j.$
 From (3 .9) we see that
  \begin{equation}
  R^{\ j} _{i \ kl}=-R^{\ j} _{i \ lk}\quad and \quad Q^{\ j} _{i \
  lk}=-Q^{\ j} _{i \ kl}.
  \end{equation}
If $D$ is a torsion-free, then $Q=0$.  Differentiating (3 .2), then
we have the first Bianchi identity
   \begin{equation}
  \omega^i \wedge \Omega^j_i=0,
\end{equation}
  which implies the first {\it Bianchi identity }for $R$:
 \begin{equation}
         R^{\ j} _{i \ kl} +R^{\ j} _{k \ li}+R^{\ j} _{l \ ik}=0,
    \end{equation}
    and
       \begin{equation}
    P^{\ j} _{i \ kl}=P^{\ j} _{k \ il}.
    \end{equation}
    The Exterior differentiation of (3 .8) gives rise to the {\it Second Bianchi
    identity}:
    \begin{equation}
    d\Omega^{\ j}_ i -\omega ^{\ k}_i \wedge \Omega^{\ j}_ k  +\omega
^{\ j}_k \wedge\Omega^{\ k}_
    i=0.
    \end{equation}
    We decompose the covariant derivatives of the Cartan tensor on
     $TM$
  \begin{equation}
    dA_{ijk}-A_{ljk} \omega ^{\ l} _i - A_{ilk}\omega ^{\ l} _j
-A_{ijl}\omega ^{\
    l} _k =A_{ijk | l}\omega ^l +A_{ijk.l}\omega ^{n+l},
    \end{equation}
    and in the similar way $\forall m\in\mathbb{N}$, for
    $\overset{_m}{A}_{ijk}$
    we have:
    \begin{equation}
    d \overset{_m}{A}_{ijk}-\overset{_m}{A}_{ljk}\omega ^{\ l}
     _i-\overset{_m}{A}_{ilk}
    \omega
^{\ l} _j-\overset{_m}{A}_{ijl}\omega ^{\ l}_k
=\overset{_m}{A}_{ijk|l}\omega ^l +\overset{_m}{A}_{ijk.l}\omega
^{n+l}.
    \end{equation}
 Clearly from (3 .15) and (3 .16), we find that for each {\it l} and
 $ \forall m\in\mathbb{N}$
 \begin{eqnarray}
A_{ijk |l}, A_{ijk .l},\ \overset{_m}{A}_{ijk|l}\quad \textrm{and}
\quad \overset{_m}{A}_{ijk.l},
\end{eqnarray}
are symmetric in i, j, k. Put
$\overset{_m}{A}_{ijk}=\overset{_m}{A}( e_i,e_j,e_k)$ and
$\overset{_m}{A^l}_{ij}=g^{kl} \overset{_m}{A}_{ijk} $, $\forall
m\in\mathbb{N}$ . By definition of
 $\dot A$ and $\overset{_m}{A}$ one has,
  \begin{equation}
  A_{ijk |n}= \dot A_{ijk},
 \end{equation}
where we use the notation $\overset{_m}{A}_{ijk
|n}=\overset{_m}{A}_{ijk |s}\ell^s$
  for all $m \in \mathbb{N}$ and
 \begin{equation}
   \overset{_m}{A}_{ijk|n}= \overset{_{m+1}}{A}_{ijk}.
 \end{equation}
 It follows from (3 .15)
 \begin{equation}
   A_{njk | l}=0 \quad , \quad  A_{njk .
   l}=-A_{jkl},
     \end{equation}
     and from (3 .16) we have
\begin{equation}
  \forall m\in\mathbb{N},\quad \overset{_m}{A}_{njk |l}=0 \quad , \quad
\overset{_m}{A}_{njk .l}=-\overset{_m}{A}_{jkl}.
     \end{equation}
 In this relation the following results
are well known: \setcounter{alphthm}{0}
 \begin{propriete}{\emph{(\cite{Ca}, \cite{M})}} Let $(M,F)$ be a Finsler manifold.
  Then for the Cartan connection (or Hashiguchi
connection), hv-curvature $P^{\ i} _{j \ kl}=0$ if and only if F is
a Landsberg metric.
\end{propriete}
\begin{propriete}{\emph{(\cite{BCS})}}Let $(M,F)$ be a Finsler manifold.
Then for the Chern connection (or Berwald connection), hv-curvature
$P^{\ i} _{j \ kl}=0$ if and only if F is a Berwald metric.
\end{propriete}
\begin{propriete}
{\emph{(\cite{Sh1})}} Let $(M,F)$ be a Finsler manifold. Then for
the Shen connection, hv-curvature $P^{\ i} _{j \ kl}=0$ if and only
if F is Riemannian.
\end{propriete}
Analogously  we have the following result.
\begin{thm} Let $(M,F)$ be a Finsler manifold. Then for the  Berwald-type
 connection, hv-curvature $P^{\ i} _{j \ kl}=0$ if and only if F is a
 Berwald metric.
\end{thm}
\vskip 0.4 true cm

\begin{proof} Let $(M,F)$ be a Finsler manifold. Differentiating
(3 .3),
 and using (3 .2), (3 .3), (3 .8), (3 .15), (3 .16), (3 .17), (3 .18)
 , (3 .19), (3 .20) and (3 .21)
  leads to
 \begin{eqnarray}
g_{kj} \Omega^{\ k}_i+g_{ik}\Omega^{\
k}_j&=&\nonumber-2A_{ijk}\Omega^k_n -2A_{ijk | l} \omega ^l \wedge
\omega^{n+k}+2A_{ijk.l} \omega ^{n+k} \ \wedge \omega ^{n+l}\\&+&
\nonumber k_{_{1}}(\dot A_{ijk | l}\omega ^l+\dot A_{ijk . l}\omega
^{n+l})\wedge \omega ^k+\cdots
\\ &+&k_{_{m}}(\overset{_m}{A}_{ijk |l}\omega ^l+\overset{_m}{A}_{ijk .l}\omega
^{n+l})\wedge \omega ^k.
  \end{eqnarray}
By using (3 .9) and (3 .22) we gives the following
 \begin{eqnarray}
 \nonumber R_{ijkl}+R_{jikl}&=&2k_{_1}\left\{\dot A_{ijl |k}-
 \dot A_{ijk | l}\right\}+\cdots+2k_{_m}\left\{ \overset{_m}{A}_{ijl |k}-
 \overset{_m}{A}_{ijk |l}\right\}\\&-&2A_{ijs}R^{\ s}_{n \ kl},
  \end{eqnarray}
   \begin{equation}
   P_{ijkl}+P_{jikl}=-2\{k_{_1}\dot A_{ijk.l}+\cdots+k_{_m}
  \overset{_m}{A}_{ijk .l}\}
   -2A_{ijl| k}-2A_{ijs}P^{\ s}_{n \ kl},
    \end{equation}
  \begin{equation}
A_{ijk.l}=A_{ijl.k}.
   \end{equation}
 Permuting $i,j,k$ in (3 .24) yields
     \begin{eqnarray} P_{ijkl}&=&-\{k_{_1}\dot A_{ijk.l}+\cdots+
     k_{_m}\overset{_m}{A}_{ijk .l}\}
     -(A_{ijl | k}+A_{jkl | i}-A_{kil | j})\nn \\&&+A_{kis}
P^{\ s} _{n \ jl}-A_{jks} P^{\ s} _{n \ il}-A_{ijs} P^{\ s} _{n\
kl},
  \end{eqnarray} and
\begin{equation}
P_{njkl}=\{k_{_1}\dot A_{jkl}+\cdots+
k_{_m}\overset{_m}{A}_{jkl}\}-\dot A_{jkl},
\end{equation}
that is because $P_{njnl}=0.$ Now if F is Berwald metric from (3
.26) and (3 .27)
 we conclude $P=0$.\\
 Conversely let $P=0$. It follows
from (3 .27),
\begin{equation}
k_{_1}\dot A_{jkl}+\cdots+k_{_m}\overset{_m}{A}_{jkl}=\dot A_{jkl},
\end{equation}
  By means of (3 .26) we have
$$
k_{_1}\dot A_{ijk.l}+\cdots+k_{_m}\overset{_m}{A}_{ijk .l}=A_{kil |
j}-A_{ijl | k}-A_{jkl | i}.$$
 Permuting $i, j, k$ in the above identity yields
$$k_{_1}\dot A_{ijk.l}+\cdots+k_{_m}\overset{_m}{A}_{ijk .l}=A_{jkl | i}-A_{kil | j}
-A_{ijl | k},$$ then
$$
A_{ijl | k}=A_{jkl | i}.$$ Letting $k=n$ in the above relation,  we
can conclude \begin{equation}
 \dot A_{ijk}=0.
 \end{equation}
 It is obvious that
\begin{equation}
\forall m\in\mathbb{N},\quad \overset{_m}{A}_{ijk}=0.
 \end{equation}
 Therefore from (3 .24), (3 .26), (3 .27) and (3 .30) we conclude
that $A_{ijk |l}=0$, thus $F$ is
 Berwald metric.\end{proof}
\section{\bf {\bf \em{\bf Some Applications}}}
\subsection{ Preliminaries on geodesics and completeness.}
 In this section we explore the notion of  geodesics
to introduce the  concept of completeness for Finsler manifolds.
 Let $c:[a,b] \rightarrow M$ be a
unit speed $C^\infty$ curve in (M,F). The canonical lift of $c$ to
$TM_0$ is defined by
\[
\hat c:={dc \over{dt}}\in TM_0.
\]
It is easy to see that $
 \rho({d{ \hat c} \over {dt}})=\ell_{\hat c},
$ where $c$ is called a {\it geodesic} if its canonical lift $\hat
c$ satisfies
\[
 {d{ \hat c} \over {dt}}= {\overline{\ell}_{\hat c}},
\]
where $\bar \ell$ is the geodesic field on
$TM_0$ defined for $\ell \in HTM$ by $\rho {(\bar \ell)}=\ell$. \\
 Let $I_xM=\{v \in T_xM,F(v)=1\}$ and
$IM=\bigcup_{p \in M}I_xM$.
 Where $I_xM$ is called  the {\it indicatrix}, and it is
a compact set.
 We can show that the projection of integral curve $ \varphi(t)$ of
$\bar \ell$ with $ \varphi(0)\in IM$ is a unit
speed geodesic $c$ whose canonical lift is $\hat c(t)= \varphi(t)$.\\

A Finsler manifold $(M,F)$ is said to be {\it backward geodesically
complete} (or {\it forward geodesically complete}) if every geodesic
$c(t)$, $a\leq t< b$ ($a< t\leq b$), parameterized to have constant
Finslerian speed, can be extended to a geodesic defined on $a\leq t<
\infty$ ($-\infty< t\leq b$).  A Finsler manifold (M,F) is said to
be {\it complete} if it is both forward and backward geodesically
complete. \\
   Let $c$ be a unit speed geodesic in M. A section $X=X(t)$
of  $\pi ^* TM$ along $\hat c$ is said to be parallel if $D_{d{ \hat
c} \over {dt}}X=0$. For $v \in TM_0$, define $\| A \|_v=sup
A(X,Y,Z)$ and $\|\dot A \|_v=sup \dot A(X,Y,Z)$, where the supremum
is taken over all unit vectors of $ \pi_v ^* TM$. Put $\| A
\|_v=sup_{v \in IM}\| A\|_v$ and $\|\dot A
\|_v=sup_{v \in IM}\|\dot A\|_v$ .\\
\subsection{Application of Berwald-type connections.}
In this subsection we are going to use two especial cases of
Berwald-type connections introduced in section 2.
 A useful property of this connection is, its adaptive form for applying to
  the different applications. In fact one can use a suitable special case
of this connection to find a geometric interpretation for solutions
of some differential equations formed by Cartan tensor and its
derivatives in Finsler spaces. For example we prove the following
theorem.
\begin{thm}  Let (M,F) be a complete Finsler manifold with bounded
Landsberg tensor. Then $F$ is a Landsberg  metric if and only if
$P_{jkl}=0$.
\end{thm}
\vskip 0.4 true cm

\begin{proof}  To prove this theorem we introduce a connection
 for which we have
  put $k_{_1}=k_{_3}=\cdots=k_{_m}=0$ and $k_{_2}\neq0$ in (3.27).
Let F be a Landsberg metric,
  then from (3.27) we find
that $P_{jkl}=0$. Conversely if $P_{jkl}=0$ then, we have following
differential equation:
\begin{equation}
k_{_{m}}A^{^{(m)}}+\cdots+k_{_{2}}A^{^{(2)}}+(k_{_{1}}-1)\dot A=0.
\end{equation}
If $k_{_1}=k_{_3}=\cdots=k_{_m}=0$ and $k_{_2}\neq0$ then we find an
special Berwald-type connection for which we have
\begin{equation}
k_{_{2}}\ddot A-\dot A=0.
\end{equation}
On the other hand
\begin{equation}
\frac{d \dot A}{dt}=\ddot A.
\end{equation}
 We have $\dot A=e^{k_{_{\texttt 2}}t}\dot A(0)$. Using
 $\|\dot A\| < \infty$
 , and letting  $t \rightarrow  +\infty$,
then $ \dot A(0)=\dot A(X,Y,Z)=0$,
 or $ \dot A=0$ i.e., $F$ is a Landsberg metric.
 \end{proof}
By mean of the Theorem 3,  every compact Finsler manifold is a
Landsberg space if and only if $P_{jkl}$ vanishes. Next we consider
a special  Berwald-type connection and give another proof for the
following well-known result. \setcounter{thm}{0}
\begin{cor}
Let (M,F) be a complete Finsler manifold with negative constant flag
curvature and bounded Cartan tensor. Then $F$ is  Riemannian.
\end{cor}
\vskip 0.4 true cm

\begin{proof} Let $(M,F)$ be a complete Finsler manifold with
constant flag curvature $\lambda$. If $\lambda\neq 0$ we put in
(3.27) $k_{_2}=k_{_4}=\cdots=k_{_m}=0$ \ ,  $ k_{_1}=2$ and
$k_{_3}=\frac{1}{\lambda} \neq0$. We obtain a connection for which
the hv-curvature $P$ become

 \begin{eqnarray} P_{ijkl}&=&-\{2\dot A_{ijk.l}+
     \frac{1}{\lambda}\dddot A_{ijk.l}\}
     -(A_{ijl | k}+A_{jkl | i}-A_{kil | j})\nn \\&&+A_{kis}
P^{\ s} _{n \ jl}-A_{jks} P^{\ s} _{n \ il}-A_{ijs} P^{\ s} _{n\
kl},
  \end{eqnarray} and
\begin{equation}
P_{njkl}=\frac{1}{\lambda} \ \dddot A +\dot A.
\end{equation}
As $M$ has constant flag curvature we have \begin{equation}\ddot
A+\lambda A=0.\end{equation} From which we have
$P_{njkl}=\frac{1}{\lambda} \ \dddot A +\dot A=0$. By solving this
differential equation we find
\begin{equation}
A(t)=c_1+c_2 e^{\sqrt{-\lambda}t}+c_3 e^{-\sqrt{-\lambda}t}.
\end{equation}
By the assumption that the Cartan tensor is bounded, and letting
$t\rightarrow \infty$ and $t\rightarrow -\infty$, we see that
$c_2=c_3=0$. Then $A=c_1$ therefore $\dot A=0$ and F is a Landsberg
metric. From (4 .6), it is easy to see that $A=0$.
\end{proof}

\vskip 0.4 true cm

\begin{center}{\textbf{Acknowledgments}}
\end{center}
The authors should express their sincere  gratitude to Professor Z. Shen for his remarks and encouragements.\\ \\


\bigskip
\bigskip

{\footnotesize \pn{\bf Behroz Bidabad}\; \\ {Department of
Mathematics and Computer Sciences}, {Amirkabir University of Technology (Tehran Polytechnic),} {15914, Tehran, 15914 Iran.}\\
{\tt Email: bidabad@aut.ac.ir}\\

\footnotesize \pn{\bf Akbar Tayebi}; \\ {Department of Mathematics }, {Qom University,} {Qom, Iran.}\\
{\tt Email: akbar.tayebi@gmail.com}\\


\begin{thebibliography}{20}
\bibitem{Ak1} H. Akbar-Zadeh, {\it Sur Les Espaces De Finsler A Courbures Sectionnelles
Constantes}, Acad. Roy. Belg. Bull. Cl. Sci. (5), {\bf 80} (1988),
271-322.
\bibitem{Ak2}
H. Akbar-Zadeh, {\it Initiation to global Finslerian Geometry},
 North-Holland Mathematical Library, 2006.
\bibitem{BC} D. Bao and  S. S. Chern,
{\it On a notable connection in Finsler Geometry},  Houston J. of
Math. {\bf 19}(1993), 135-180.
\bibitem{BCS} D. Bao, S. S. Chern and Z. Shen,
{\it An Introduction to Riemann-Finsler Geometry}, Springer-Verlag,
2000.
\bibitem{Be}
L. Berwald, {\it Untersuchung der Kr\"{u}mmung allgemeiner
metrischer R\"{a}ume auf Grund des in ihnen herrschenden
Parallelismus}, Math. Z. {\bf 25}(1926), 40-73.

\bibitem{BT} B. Bidabad and A. Tayebi,\emph{A
 classification of some Finsler connections},  Publ. Math. Debrecen
{\bf 71}/3-4 (2007), 253-266.


\bibitem{Ca}
E. Cartan, {\it Les espaces de Finsler}, Hermann, Paris, 1934.
\bibitem{Ch}
S. S. Chern, {\it  On the Euclidean connections in a Finsler space},
Proc. National Acad. Soc., {\bf 29}(1943), 33-37; or Selected
Papers, vol. I\negthinspace I, 107-111, Springer 1989.
\bibitem{CS} S.S. Chern and Z. Shen, {\it Riemann-Finsler Geometry ,} Preprint,
2006.
\bibitem{M} M. Matsumoto, {\it Foundation of Finsler geometry and
special Finsler spaces,} Kaiseisha Press, Japan, 1986.
\bibitem{MA} R. Miron, and M. Anastasiei, {\it The Geometry of Lagrange
space: Theory and Application}, Kluwer, Dordrecht, 1994.
\bibitem{MS} M. Matsumoto and H. Shimada, {\it On Finsler spaces with
the curvarure tensors $P_{hijk}$ and $S_{hijk}$ satisfying special
conditions,} Rep. On Math. Phy. {\bf 12}(1977), 77-87.
\bibitem{Sh1} Z. Shen, {\it On a connection in Finsler Geometry},
Houston J. of Math. {\bf 20}(1994), 591-602.
\bibitem{Sh2} Z. Shen, {\it Differential Geometry of Spray and
Finsler Spaces}, Kluwer Academic Publishers, Dordrecht 2001.
\bibitem{Sh3} Z. Shen, {\it Lectures on Finsler Geometry}, Word
 Scientific, 2001.
\end{thebibliography}
\end{document}